\documentclass{amsart}

\usepackage{booktabs}
\usepackage[hidelinks]{hyperref}

\newtheorem{theorem}{Theorem}[section]

\numberwithin{equation}{section}

\newcommand{\R} {\mathbb R}
\newcommand{\M} {\mathcal M}
\newcommand{\DT} {\mathcal L}

\begin{document}

\title[Efficient computation of the Feigenbaum constants]
      {An efficient method for the computation
       of the Feigenbaum constants to high precision}

\author{Andrea Molteni}
\email{am@pass.im}

\subjclass[2010]{Primary 11Y60, 26A18}

\date{December 15, 2015}

\begin{abstract}
We propose a new practical algorithm for computing
the Feigenbaum constants $\alpha$ and $\delta$,
having significantly lower time and space complexity
than previously used methods.
The algorithm builds upon well-known linear algebra techniques,
and is easily parallelizable.
An implementation of it has been developed and used
to determine both constants to $10\,000$ decimal places.
\end{abstract}

\maketitle


\section{Introduction}\label{sec:introduction}
The Feigenbaum constants $\alpha$ and $\delta$~\cite{feigenbaum78}
arise as limits in the theory of iteration of real functions.
Their natural definition is, however, unpractical for computation to high precision,
as it leads to slowly-converging, exponential-time algorithms
(see e.g.~\cite{briggs89}).

A much more suitable characterization is introduced
by Feigenbaum himself in~\cite{feigenbaum78,feigenbaum79}.
Let us define an operator $T$, acting on functions $g:\R\rightarrow\R$, by
\begin{equation}\label{eq:operatort}
(Tg)(x) = \frac{g(g(g(1)\,x))}{g(1)}.
\end{equation}
By using a particular even analytic function $g:[-1,1]\rightarrow\R$
invariant under~$T$, having a local quadratic maximum at $x=0$
and such that $g(0)=1$, we can determine Feigenbaum's $\alpha$ constant as
\begin{equation}\label{eq:alpha1g1}
\alpha = \frac{1}{g(1)}.
\end{equation}
Feigenbaum's $\delta$ constant can then be computed
as the largest real eigenvalue of the linear operator $\DT$
defined by
\begin{equation}\label{eq:dt}
(\DT f)(x) = \alpha g'(g(x/\alpha))\cdot f(x/\alpha) + \alpha f(g(x/\alpha)).
\end{equation}
It is therefore apparent that, within the framework that we have just outlined,
an efficient determination of the fixed-point function $g$
is key to a precise estimation of both $\alpha$ and $\delta$.

A description of our method for computing the Feigenbaum constants
is given in Section~\ref{sec:method} of this paper.
In Section~\ref{sec:analysis} we analyze the complexity of the method
while comparing it to other methods that have appeared in the literature.
Finally, Section~\ref{sec:results} discusses some implementation details
and reports a few numerical results.


\section{The new method}\label{sec:method}
Let $n$ be an integer greater than $1$.
In order to approximate the function $g$, fixed point of~\eqref{eq:operatort},
we model it as a truncated Chebyshev series of the form
\begin{equation}\label{eq:gcheby}
\tilde g_n(x) = \sideset{}{'} \sum_{j=0}^{n-1} c_j T_{2j}(x),
\end{equation}
where the prime indicates the standard convention that the first term of the sum
is to be halved, and where
$T_n$ is the $n$th \emph{Chebyshev polynomial of the first kind}~\cite{mason-handscomb},
defined by the recurrence relation
\begin{equation*}
\begin{aligned}
T_0(x) &= 1 \\
T_1(x) &= x \\
T_{n+1}(x) &= 2xT_n(x) - T_{n-1}(x) \quad \forall n > 0.
\end{aligned}
\end{equation*}
We intend to determine the coefficients $c_0,\ldots,c_{n-1}$ by a collocation method.
To this end, we require that $\tilde g_n$ satisfy the Feigenbaum-Cvitanovi\'c equation
\begin{equation}\label{eq:feig-cvit}
g(1)\,g(x) - g(g(g(1)\,x)) = 0,
\end{equation}
which characterizes fixed points of~\eqref{eq:operatort},
and evaluate the resulting equation at the \emph{Chebyshev nodes}
\begin{equation}\label{eq:chebnodes}
t_i = \cos\left(\frac{2i-1}{4n}\,\pi\right), \quad i=1,\ldots,n,
\end{equation}
obtaining the system of $n$ nonlinear equations
\begin{equation}\label{eq:f0cheby}
F(c_1,\ldots,c_n) :=
\begin{bmatrix}
\tilde g_n(1)\,\tilde g_n(t_1) - \tilde g_n(\tilde g_n(\tilde g_n(1)\,t_1)) \\
\vdots \\
\tilde g_n(1)\,\tilde g_n(t_n) - \tilde g_n(\tilde g_n(\tilde g_n(1)\,t_n))
\end{bmatrix}
=
\begin{bmatrix}0 \\ \vdots \\ 0\end{bmatrix}.
\end{equation}

As an initial approximation of the solution of~\eqref{eq:f0cheby} we take
the $n$-tuple $x^{(0)}=(\hat c_0,\ldots,\hat c_{n-1})$, where the first
\begin{equation}\label{eq:msqrt}
m \approx \frac32 \sqrt{n}
\end{equation}
coefficients $\hat c_0,\ldots,\hat c_{m-1}$ have been obtained
by applying the present method to determine $\tilde g_m$ with a collocation of $m$ points,
and where we let $\hat c_j = 0$ for $j=m,\ldots,n-1$.
For the case $m=2$, we use $(\hat c_0,\hat c_1)=(0.6,-0.7)$.

By finite differences, we compute the Jacobian matrix $B_0$ of $F$ at $x^{(0)}$
to the same precision of the coefficients $\hat c_j$;
we then explicitly invert this matrix by Gaussian elimination
to obtain $B'_0 = B_0^{-1}$.
Finally, we apply
the Inverse Column Updating Method (ICUM) by Mart\'inez and Zambaldi~\cite{icum}
to iteratively solve equation~\eqref{eq:f0cheby}.
For $k\geq 0$ we compute
\begin{equation}
x^{(k+1)} = x^{(k)} - B'_k F(x^{(k)}),
\end{equation}
and define
\begin{align*}
s^{(k)} &= x^{(k+1)} - x^{(k)}, \\
y^{(k)} &= F(x^{(k+1)}) - F(x^{(k)}).
\end{align*}
We then choose $j_k\in\{1,\ldots,n\}$ such that
\begin{equation}
|y^{(k)}_{j_k}| = \|y^{(k)}\|_\infty := \max\bigl\{|y^{(k)}_1|,\ldots,|y^{(k)}_n|\bigr\},
\end{equation}
and update our approximate inverse Jacobian by the formula
\begin{equation}
B'_{k+1} = B'_k + \frac{s^{(k)} - B'_k y^{(k)}}{y^{(k)}_{j_k}} \otimes e^T_{j_k},
\end{equation}
where $e_{j_k}$ denotes the $j_k$-th element of the canonical basis of~$\R^n$.

As soon as $\tilde g_n$ has been determined to sufficient precision,
an approximation of $\alpha$ can be easily computed as $1/\tilde g_n(1)$.


To determine $\delta$, we apply the Arnoldi iteration~\cite{arnoldi}
to an $n\times n$ real matrix $L$ approximating
the infinite-dimensional operator $\DT$ defined by~\eqref{eq:dt}.
However, we never explicitly construct $L$;
instead, we exploit the ``black-box'' nature of the Arnoldi process,
which only requires the ability to determine the product $Lv$ for any vector $v\in\R^n$,
and has no need to access or manipulate the entries of $L$ directly.
For any given $v\in\R^n$, whose entries represent the Chebyshev coefficients
of an even analytic function $f:[-1,1]\rightarrow\R$,
\begin{equation*}
f(x) = \sideset{}{'} \sum_{j=0}^{n-1} v_j T_{2j}(x),
\end{equation*}
we compute $Lv$
by evaluating~\eqref{eq:dt} at the $n$ Chebyshev nodes~\eqref{eq:chebnodes},
and by applying a discrete cosine transform
to infer the Chebyshev coefficients of~$\DT(f)$.

The Arnoldi iteration consists in a Gram-Schimdt-like process for
reducing a generic matrix to Hessenberg form
(where all entries below the first subdiagonal are zero)
while preserving its spectrum.
It belongs to a class of linear algebra algorithms,
known as \emph{iterative methods},
that provide a meaningful partial result at each iteration,
in contrast to \emph{direct methods}, which only give a useful result upon completion.
In particular, the $k$th iteration of the Arnoldi method
yields a $k\times k$ Hessenberg matrix $H_k$,
whose eigenvalues are known as the \emph{Ritz eigenvalues}.
It is often observed in practice that these eigenvalues converge
to the extreme eigenvalues of the input matrix as $k$ tends to $n$.
In our case we found that the extreme eigenvalue of $H_k$ converges
to the extreme eigenvalue of $L$ to maximum precision after only
\begin{equation}\label{eq:p}
p \approx {3\sqrt{n}}
\end{equation}
iterations
if we choose as starting vector for the Arnoldi process $e_1$,
the first element of the canonical basis of~$\R^n$.
In fact, this vector constitutes a good enough approximation of an eigenvector
relative to the extreme eigenvalue of $L$,
whose actual entries are observed to be exponentially decreasing in absolute value.

At each iteration, $\delta$ can be quickly estimated as
a root of the characteristic polynomial of $H_k$, $p_k(t) = \det(tI_k-H_k)$,
by applying the classical secant method.
This calculation is numerically stable,
as $p_k(t)$ can be evaluated directly by exploiting the Hessenberg structure of $H_k$,
without ever explicitly computing the coefficients of the polynomial.


\section{Analysis of the method}\label{sec:analysis}
The classical approach~\cite{feigenbaum78,feigenbaum79,briggs91}
to the numerical approximation of the function $g$
models it as a truncated power series of the form
\begin{equation}\label{eq:gpower}
\bar g_n(x) = 1 + \sum_{j=1}^n a_jx^{2j}.
\end{equation}
A set of $n$ nonlinear equations in the $n$ unknowns $a_1,\ldots,a_n$
similar to~\eqref{eq:f0cheby}
is obtained by a collocation method
similar to the one described in Section~\ref{sec:method},
and solved by applying an $n$-dimensional Newton's method.
The higher $n$, the higher the accuracy of the approximation of $g$.
In particular, it is observed that the number of correct digits
of the resulting approximation of $\alpha$ increases about linearly with $n$.
Contrary to what suggested in~\cite{briggs91}, we found that in practice
$n$ need never be greater than the number of decimal digits desired for $\alpha$.

What we have just summarized is the classical method referenced
in most of the literature we have consulted,
and that we will consider as a baseline to compare our proposed improvements against.
To do so, we first need to analyze the performance of this method.

In what follows, let $\M(n)$ denote a function such that
two numbers of length $n$ can be multiplied in time $O(\M(n))$.
For example, $\M_{S}(n) = n\log n\log\log n$ can be one such function
if using the Sch\"onhage-Strassen algorithm~\cite{schonhage-strassen}.
It is well-known that addition, subtraction and division of numbers of length $n$
can also be performed in time $O(\M(n))$.
Since we are interested in calculations to thousands of decimal places,
we can assume
\begin{equation}\label{eq:mulbounds}
O(n) < O(\M(n)) \leq O(n^{\log_2 3}),
\end{equation}
where the upper bound is given by Karatsuba's algorithm~\cite{karatsuba}.

The execution time of one iteration of Newton's method as proposed above is dominated
by the computation of an $n\times n$ Jacobian matrix
and by the solution of a linear system of $n$ equations in $n$ unknowns.
Since there does not seem to be a cheap way to evaluate the partial derivatives
analitically, the entries of the Jacobian matrix
are approximated by finite differences. This requires for every entry
at least three evaluations of~\eqref{eq:gpower},
which can be accomplished in $O(n\,\M(n))$ by using Horner's method.
Therefore, computing the whole Jacobian matrix at a given set of coefficients
$\{a_j\}$ takes time $O(n^3\M(n))$. The same time requirement applies to
solving the aforementioned $n\times n$ linear system, as is well known.

If the initial approximation of the function $\bar g_n$ is chosen close enough
to the actual solution, Newton's method converges quadratically.
The total time required to approximate $g$ by $\bar g_n$
to maximum precision is therefore
$O(n^3\log n\,\M(n))$, while the memory requirements are $O(n^3)$,
dominated by the storage costs of the Jacobian matrix.

The fundamental issue with modeling $g$ as in~\eqref{eq:gpower}
is that it leads to Jacobian matrices that are very ill-conditioned
near the solution of the nonlinear system.
While this does not prevent root-finding methods from converging,
it imposes strong requirements on the precision to which the Jacobian
matrix has to be computed, increasing the computational burden.
This led us to abandon model~\eqref{eq:gpower} in favor of model~\eqref{eq:gcheby},
as it is well known from interpolation theory that
the Chebyshev polynomials $T_n$ form a much more numerically stable basis
for the space of polynomials than the basis given by monomials.
And while the sum in~\eqref{eq:gcheby} may look more complex
than that in~\eqref{eq:gpower},
it can still be evaluated in time $O(n\,\M(n))$
by using Clenshaw's algorithm~\cite{mason-handscomb,nr}.

An immediate consequence of this change of basis
is a uniform increase in the accuracy of the resulting approximations.
In particular, we observed an increase of about 11\%
in the number of correct digits of $\alpha$ obtained for any given $n$.

More importantly, this change drastically improves the conditioning
of the Jacobian matrix, which in turn allows us to relax
the precision requirements on its approximation.
Experimentally we found that,
if we use a quasi-Newton method like Broyden's~\cite{broyden}
to solve the nonlinear system of equations,
it is often possible to reduce the precision
of the approximate Jacobian to $O(\sqrt{n})$
while retaining linear convergence.
This reduces the time complexity of computing the initial
approximate inverse Jacobian to $O(n^3\M(n^{0.5}))$,
and the total memory requirements to $O(n^{2.5})$.

We have investigated several quasi-Newton methods~\cite{broyden,projup,thomas,colum,icum,sectens}
for solving~\eqref{eq:f0cheby},
and in practice we found ICUM to be particularly effective,
as it requires only one matrix-vector product per iteration
and its update formula is significantly cheaper than that of other
multidimensional secant methods such as Broyden's.
If we start from the $m$-point approximation
described in Section~\ref{sec:method},
ICUM reaches convergence in~$O(n/m)=O(\sqrt{n})$ iterations.
Each iteration is dominated in time by the evaluation of $F$,
which involves $3n$ evaluations of~\eqref{eq:gcheby},
i.e.\ $O(n^2)$ basic operations at full precision,
so the total time complexity of the root-finding method
amounts to $O(n^{2.5}\M(n))$.

We deem it important to note that
the explicit inversion of $B_0$ appears to be the best option in this case.
In our experiments with quasi-Newton methods,
alternatives such as QR decomposition
increased computational and/or storage costs
without bringing any significant benefit.

Finally, we observe that,
since only $m$ of the $n$ coefficients of the initial approximation of~\eqref{eq:gcheby}
are nonzero, the actual time needed to compute $B_0$
by finite differences can be lowered to $O(n^{2.5}\M(n^{0.5}))$
by factoring the evaluation of isolated trailing terms in~\eqref{eq:gcheby}
out of Clenshaw's algorithm.
This is because an isolated $T_{2j}(x)$ can be computed with just
$O(\log j)$ multiplications and additions by recursively exploiting the relations
\begin{equation}
\left\{
\begin{aligned}
T_{2k}(x) &= 2\,T_k(x)\,T_k(x) - T_0(x) \\
T_{2k-1}(x) &= 2\,T_k(x)\,T_{k-1}(x) - T_1(x)
\end{aligned}
\right.
\qquad \forall k > 0,
\end{equation}
which follow from the basic property of Chebyshev polynomials
\begin{equation*}
2\,T_m(x)\,T_n(x) = T_{m+n}(x) + T_{|m-n|}(x) \qquad \forall m,n.
\end{equation*}

\begin{table}[t]
\centering
\caption{Complexity of the new method for approximating $g$.}
\label{tab:gsteps}
\begin{tabular}{lp{80mm}}
\toprule
Time & Step \\
\midrule
$O(n^{1.25}\M(n^{0.5}))$ & Approximate $g$ on a collocation of $m$ nodes \\
$O(n^{2.5}\M(n^{0.5}))$ & Compute the initial approximate Jacobian~$B_0$
to~reduced precision by~finite differences \\
$O(n^3\M(n^{0.5}))$ & Invert $B_0$ in place by Gaussian elimination \\
$O(n^{2.5}\M(n))$ & Apply a quasi-Newton root-finding method \\
\bottomrule
\end{tabular}
\end{table}

Table~\ref{tab:gsteps} summarizes the steps of the method we developed
for solving~\eqref{eq:f0cheby}.
Inequality~\eqref{eq:mulbounds} entails that
the last step, namely the multidimensional secant method,
dominates the other steps timewise.


We will now assume having computed a high-precision approximation of $g$
and proceed to analyze methods for computing $\delta$ as an eigenvalue
of the linear operator $\DT$ defined in~\eqref{eq:dt}.

The classical approach appearing in the literature involves
approximating the infinite-dimensional linear eigenvalue problem
by explicitly constructing an $n\times n$ matrix $L$ and
applying standard linear algebra algorithms for the determination of eigenvalues.
For example, the power method~\cite{golub-vanloan}
is mentioned by Briggs in~\cite{briggs91}.
The matrix $L$ can be computed in many different bases
for the space of even polynomials.
Traditionally monomials have been used~\cite{briggs91},
but we prefer Chebyshev polynomials
as they provide better stability without introducing any significant drawback.

It is easy to see that, in both cases,
the total time requirement for building $L$ is $O(n^3\M(n))$.
Each iteration of the power method calls then for a matrix-vector product,
and the convergence of the method is only linear,
for a total time complexity of $O(n^3\M(n))$.

As explained in Section~\ref{sec:method},
our method avoids the explicit formation of $L$
and instead considers its action on even functions
represented by the $n$ coefficients of their truncated Chebyshev series expansion.
For any such function,
the evaluation of~\eqref{eq:dt} at the $n$ Chebyshev nodes takes time $O(n^2\M(n))$;
the Chebyshev coefficients of the transformed function
can then be retrieved by a fast cosine transform~\cite{nr}
in time $O(n\log n\,\M(n))$.

The $k$th iteration of the Arnoldi method,
yielding a $k\times k$ Hessenberg matrix $H_k$,
involves one evaluation of $\DT$ and $O(k)$ inner products of vectors of $\R^n$,
for a total time requirement of $O((n^2+kn)\,\M(n))$.

At each iteration, the characteristic polynomial of $H_k$
can be evaluated with just $O(k^2)$ elementary operations
by exploiting the Hessenberg structure of $H_k$.
Thus, by applying the classical, superlinearly convergent secant method
to this polynomial
we can produce an estimate of $\delta$ in time $O(k^2\log n\,\M(n))$.
For $k\ll n$, and in particular for $k\leq p$,
this computation is much cheaper than one iteration of the Arnoldi method,
and can therefore be performed after each iteration
without incurring any significant penalty.

Since the approximation of $\delta$ reaches the maximum precision
allowed by model \eqref{eq:gcheby} after just $p\approx 3\sqrt{n}$ iterations,
we conclude that our method completes in time $O(n^{2.5}\M(n))$.
The memory requirements are dominated by the $p$ Ritz vectors
maintained by the Arnoldi iteration,
so the total space complexity is, once again, $O(n^{2.5})$.


\section{Numerical results}\label{sec:results}
The methods presented in this paper rely on just a few
relatively time-consuming algorithms.
Most of the linear algebra algorithms we mentioned,
such as matrix-vector multiplication and Gaussian elimination,
are easily parallelizable.
The same holds for the discrete cosine transform, and
the approximation of a Jacobian matrix by finite differences
is even embarrassingly parallel.

One problem arises, however, when trying to parallelize the Arnoldi iteration.
In its standard implementation, the inner loop of this algorithm
consists of a modified Gram-Schmidt process.
This is often required to avoid the numerical instability
intrinsic to the classical Gram-Schmidt process,
but in a parallel setting it also introduces the need for
frequent communication of large amounts of data between nodes.
It is worth noting that, in our case, stability did not prove to be an issue,
so this problem could simply be fixed by using the classical Gram-Schmidt algorithm.

We will now report some numerical results.
For comparison,
the best previous estimate of the Feigenbaum constants known to us
is that of Broadhurst~\cite{broadhurst99},
who computed both $\alpha$ and $\delta$ to 1018 decimal places
by using a collocation of 700 points and 400~MB of memory.
In 1999, this calculation took 3 days.

A parallel implementation of the method described
in Section~\ref{sec:method}
was developed in C using the GNU MPFR library~\cite{mpfr} and MPI~\cite{mpi}.
We were able to replicate Broadhurst's results with a collocation
of 630 Chebyshev nodes, and using only 33~MB of memory.
This calculation took less than one minute on a modern desktop computer;
on the same machine, an implementation of the more classical methods
took over one hour.

We used our methods to compute estimates of $\alpha$ and $\delta$
for all values of $n$ from $2$ to $1000$.
Let us denote by $\tilde\alpha_n$ and $\tilde\delta_n$ the approximations
obtained on $n$ nodes for $\alpha$ and $\delta$ respectively.
The quantities
\begin{align}
D_\alpha(n) &= \log_{10}\frac{|\alpha|}{|\tilde\alpha_n-\alpha|}, &
D_\delta(n) &= \log_{10}\frac{|\delta|}{|\tilde\delta_n-\delta|}
\end{align}
represent the number of correct decimal digits of these approximations,
and constitute a measure of their accuracy.
As one may expect, $D_\alpha(n) > D_\delta(n)$ holds for almost all $n$,
since the computation of $\tilde\delta_n$ involves a more elaborate process
than the computation of $\tilde\alpha_n$,
and both are based on the same approximation of $g$.
However, the difference $D_\alpha(n) - D_\delta(n)$ was in all cases very small,
never exceeding~$4$.
Based on the data we collected, we conjecture that this difference
grows proportionally to~$\log n$.

It is also observed that
\begin{equation}
\liminf_{n\rightarrow\infty} \frac{D_\delta(n)}{n} \approx 1.63,
\end{equation}
converging from below.
This means that the number of correct digits in our approximations
grows at least linearly with~$n$.
When using the monomial basis, the same limit is about~$1.46$.

Our best estimate of the Feigenbaum constants was obtained
with a collocation of $6144$ points, and checked with $6160$ points.
Each run took 3 days using 7~GB of RAM of a modern desktop computer.
The results matched to $10\,026$ decimal places for $\alpha$
and $10\,022$ decimal places for $\delta$, respectively.

We used these estimates to test the transcendence of $\alpha$ and $\delta$
with the PSLQ algorithm by Ferguson and Bailey~\cite{pslq}.
This is an integer relation algorithm that can be used
to determine whether a given real number $x$ is likely to be algebraic
by searching for integers $m_0,\ldots,m_n$, not all zero, such that
$\sum_{i=0}^n m_ix^i = 0$.
In case that no such relation is found, the algorithm provides a lower bound
on the norm of any potential tuple $(m_0,\ldots,m_n)$ satisfying the equation.
We obtained the following result.

\begin{theorem}\label{th:pslq}
If either $\alpha$ or $\delta$ is a root of an integer polynomial
of degree~$120$ or less,
then the Euclidean norm of the coefficients exceeds~$1.5\times 10^{79}$.
\end{theorem}

Theorem~\ref{th:pslq} constitutes a significant improvement
on the bounds published by Briggs in~\cite{briggs97}
for polynomials of degree up to $20$.


\section{Conclusion}\label{sec:conclusion}
We introduced a new method for the precise approximation of
the universal Feigenbaum function $g$
and of the Feigenbaum constants $\alpha$ and $\delta$.
This method reduces time complexity from $O(n^3\log n\,\M(n))$ to $O(n^{2.5}\M(n))$,
and space complexity from $O(n^3)$ to $O(n^{2.5})$.


\bibliographystyle{amsplain}
\bibliography{hifeig}

\end{document}